\documentclass{article}
\usepackage{color,times,amsmath,amsfonts,latexsym,epsfig,epsf,colordvi}
\usepackage{url,hyperref}
\usepackage[english]{babel}

\usepackage{xypic}

 \usepackage{tikz}
\usetikzlibrary{arrows}
\usepackage{graphicx}
\usepackage{amssymb}
\usepackage{psfrag}
\usepackage{subfigure}
\usepackage{float}
\usepackage{comment}

\usepackage{lipsum}
\usepackage{epstopdf}
\usepackage{algorithmic}
\usepackage{amsopn}
\usepackage{amsfonts}
\usepackage{multirow}
\usepackage{mathptmx}
\usepackage{amsmath}
\usepackage{psfrag,epsfig}
\usepackage{graphicx}
\usepackage{subfigure}
\usepackage{CJK}

\ifpdf
  \DeclareGraphicsExtensions{.eps,.pdf,.png,.jpg}
\else
  \DeclareGraphicsExtensions{.eps}
\fi

\newcommand{\RR}{\mathbb {R}}

\newcommand{\bp}{\begin{pmat}}
\newcommand{\ep}{\end{pmat}}

\makeatletter
\def\adots{\mathinner{\mkern1mu\raise\p@\vbox{\kern7\p@\hbox{.}}\mkern2mu\raise4\p@\hbox{.}\mkern2mu\raise7\p@\hbox{.}\mkern1mu}}
\makeatother

\title{ \vspace*{-8mm}{ \bf The Eight Epochs of Math as regards\\  
past and future Matrix Computations \vspace*{-4mm} }}

\author{ Frank Uhlig\\
Department of Mathematics and Statistics,\\
  Auburn
University, AL 36849-5310 \ ({\tt uhligfd@auburn.edu})\vspace*{5mm}}

\begin{document}
\date{  }
\maketitle

\vspace*{-2mm}

\pagestyle{myheadings}
\thispagestyle{empty}

\vspace*{-6mm}
{\normalsize 
{\bf Abstract : }
\noindent
This talk gives a personal  assessment of Epoch making advances in Matrix Computations from antiquity and with an eye towards tomorrow.\\
 We trace the development of number systems and elementary algebra, and the uses of Gaussian Elimination methods from around 2000 BC on to current real-time Neural Network computations to solve time-varying linear equations.\\
  We include relevant advances from China from the 3rd century AD on, and from India and Persia in the 9th century and discuss the conceptual genesis of vectors and matrices in central Europe and Japan in the 14th through 17th centuries AD.\\ 
  Followed by the 150 year cul-de-sac of  polynomial root finder research for matrix eigenvalues, as well as the superbly useful matrix iterative methods and  Francis' eigenvalue Algorithm from last century.\\
  Then we explain the recent use of  initial value problem solvers to master time-varying linear and nonlinear matrix equations via Neural Networks.\\ 
We end with a short outlook upon new hardware schemes with multilevel processors that go beyond the 0-1 base 2 framework which all of our past and current electronic computers have been using.

\section*{\Large Introduction}

In this paper we try to outline the {Epoch} making achievements and transformations that have occurred over time for  computations and more specifically for matrix computations.
We will trace how  our linear algebraic  concepts and matrix computations  have progressed  from the beginning of recorded time until  today and how they will likely progress into the future. We take this limited tack simply because in modern times matrices have become the elemental and universal tools for most any computation. \\[1mm]
This  evolution of our matrix methods will be described in broad strokes.
My main emphasis is to trace the mathematical genesis of Matrices and their uses and to learn how the modern matrix  concept has evolved in the past and how it is evolving.  
I am not interested in Matrix Theory by itself, but rather in Matrix Computations, i.e.,  how matrix concepts and algorithms have been developed from approximately 3000 BC to today, and even tomorrow.\\[1mm]
This paper describes {eight noticeably separate Epochs} that are distinguished from each other  by the introduction of  evolutionary new concepts and subsequent radically new computational methods. Following the  historical trail through six historically established Epochs,  we will then look into the present and  the near future.\\[1mm]
 What drives us  to conceptualize and compute differently now, what is leading us into the 7th and possibly 8th Epoch? When and how  will we likely compute in the future. \\[1mm]
I am  not a math historian, I have never taught a class in math history. Instead throughout my academic career I have worked with matrices: in matrix theory, in applications and in numerical analysis and I  like to construct efficient new algorithms that solve matrix equations. 
The idea for this paper is in part due to my listening by chance to a very short English broadcast  from Egyptian Radio  on short wave some  40  years ago in the 1970s.
 It described an Egyptian papyrus from around { 2000 BC} that dealt with solving linear equations by   row reduction and zeroing out coefficients in systems of linear equations, i.e.,  by  what we now call ``{ Gaussian Elimination}".
  When I heard this as a young Ph. D., I was fascinated and wrote the station for more information. They never answered and when I was in Cairo many years later, the Egyptian Museum personnel could not help me either with locating the source. \\[-2mm] 

Thus I became aware that Carl Friedrich Gau\ss\ did not invent what we now call by his name; but who did? \\
For many decades, this snippet of math history just lingered in my mind until a year ago when I was sent a book  on { Neural Network} (NN) methods for solving time-varying linear and non-linear equations and  was asked to review it.
 The NN methods were - to me and my understandings of numerics then - so other-worldly and brilliant that I began to think of the incredible leaps and 'bounces' that math computations have gone through over the eons, from era to era.
 I eventually began to detect {7 or 8 computational Sea Changes}, what I call ``{ Epochs}", in our ability to compute with matrices and that is my topic.

\section{\Large A Short History of Matrix Computations}

Nobody knows how  numbers and number systems came about, just like nobody knows `who invented the wheel'. I will start  with a few historical facts about number systems and how they developed and were used across the globe in antiquity.\\[-6mm]

\subsection*{\large 1.0\quad Early Number Systems}
Humankind's first developments of number systems were very diverse and geographically widely dispersed, yet rather slow.
The first circle cypher for zero occured in Babylonia around 2500 BC, or 4500 years ago. A continent or two removed, the Mayans used the same concept and circle zero symbol  from around 40 BC. In India it was recognized during the 7th century. But zero only became recognized as {a `number to compute with'} like all the others in the 9th century in central India.
 Our  decimal system builds on the {ten numbers 0, 1, 2, 3, 4, 5, 6, 7, 8, and 9}. The decimal positional system  came from China via the Indus valley and it started to be used in Persia in the 9th century. It was combined with or  derived  from a Hindu number system of the same time period.\\[1mm]
  In fact westerners call the current decimal number symbols wrongfully   'arabic', but most westerners (and I) cannot read the license plates on cars in Egypt  since the arabic world does not use our Persian/Hindu numbers in writing but its own script using  arabic letters to designate numbers.  {Should we call our 'western' numbers 'farsi' or `hindu' instead? }\\[-4mm]
 
Various bases have been used for numbering. There have been base 2, base 8, base 10, base 12, base 16, base 60, base 200 number systems and possibly more at some time somewhere  throughout human history. 
Counting and simple computations started with {notched sticks} for record keeping and with the invention of sand or wax tablets and then the {abacus}. These simple tools were developed a little bit differently and independently in many parts of the globe.

\subsection{\large Antiquity\ \  (1st Epoch)}

Around 2200 to 1600 BC, Sumerian, Babylonian and Egyptian land survey computations became mathematized  in order to  mark and allot land after the yearly Euphrates, Tigris and Nile floods. That naturally lead to linear equations in  2, 3 or 4 variables and subsequent methods to solve them that amounted to what we now  call {row reduction} or {Gaussian elimination}.\\[1mm]
Mathematics Computations did not advance much during the Greek times as Greek mathematicians were  mainly interested in mathematical theory and in establishing the concept of a formal proof, as well as elementary number theory of which the Euclidean algorithm is still used today.\\[1mm]
 Neither did the complicated Roman numerals lend themselves to easy computations and no further computational advances happened there.
 
     \subsection{\large Early Mathematical Arts  in China, India and the Near East\ \  (2nd Epoch)}
    
 (Based in part on a lecture at Hong Kong University in 2017, given by Zhou Xiangyu, Chinese Academy of Sciences, for Chinese sources.  And on Indian and Arabic sources from elsewhere)\\[1mm]
In prehistoric and historic times (1600 BC - 1400 AD), {knot and rod calculus} were prevalent in China. They were based on a decimal positional system, the so-called "rod numerals".  These comprised the most advanced number system of the time and it was used for several millennia before being adopted and expanded in Persia and India in the 9th century AD and later on adopted in central Europe.\\[1mm]
 The {"Mathematical Classic of Sunzi"} by Suanjing (from the 3rd - 5th century) gives a detailed description of the arithmetic rules for counting rods. In the Indus valley clay tablets covered with sand were used for mathematical computations several millennia ago. Bhaskara (600 -  680 AD) in India was  first to write numbers in the Hindu positional decimal system which used the circle for zero. In 629 AD he approximated the sine function by rational expressions while commenting on Aryabhatta's (476 - 550 AD)  book "Aryabhatiyabhisya" from 499 AD. 
An Indian  contemporary of Bhaskara, Brahmagupta (598 - 665 AD) was the first to establish the rules that govern computing with zero. Brahmagupta texts were written in Sanskrit verse that used the Sanskrit word for 'eyes' to denote 2 and 'senses' for the number 5 etc. This was common in Indian mathematics and science writings at the time. The earliest record of multiplication and division algorithms using the Hindu numerals 1 through 9 and 0 was in writings by Al Khwarizmi 780 - 850 AD, a  Persian mathematician employed in Bagdad. His "Book of Manipulation and Restoration" established the golden rule of Algebra that an equation remains true if one subtracts the same quantity from both sides. He also wrote down multiplication and division rules that are identical to those of Suanjing from the 3rd to 5th century and China. To Suanjing we also owe  the Chinese Remainder Theorem. Finally the  advanced Hindu-Arabic decimal number system  was introduced into the West by Leonardo Fibonacci (1175 - 1250 AD) of Pisa in his "Liber Abaci" or 'The Book of Calculations' (1202),\\[1mm] 
 Applied and numerical computations were driving much of Chinese mathematics. Wang Xiaotang (580 - 640 AD), for example, tried to find the roots of cubic polynomials that appeared in civil engineering and water conservation problems.
In the {Mathematical Treatise in Nine Sections} of 1247, Qin Jiushao (1202 - 1261)  developed the 'Qin Jiushao method' which is now commonly called the  'Horner-Ruffini scheme' for computing with and finding roots of polynomials iteratively.
 William George Horner (1819) \cite{Horn} and Paolo Ruffini (1804-1807-1813)    reinvented the Qin Jiushao method unknowingly  six hundred years later.\\[1mm]
Chinese rod calculus was the method of choice for computing in China until the  {abacus} took over during the Ming Dynasty 1388 - 1644. Cheng Dawei 1355 - 1606 is the author of the first 'Numerical Analysis' book titled {"The General Source of Computational Methods"} and published in 1592. It describes methods to add, subtract, multiply, and divide on an   abacus. 
The abacus itself was invented in various incarnations at various times and in several locations of the globe. It essentially combines several decimal rods on one board with beads on strings. \\[1mm]
Chinese mathematicians from the 3rd century BC onwards to the 10th century AD  brought us the {"Nine Chapters on the Mathematical Arts"}  that uses the numbers 1 through 9. This book was later disseminated further to the west, to India and Persia as described above.
 In  its 7th chapter, {determinants} first appeared conceptually, while chapter 8 abstracts the concept of linear equations to  represent them by {matrix-like coefficient tableaux}.
   These {`matrix equations'} were solved in China, again  by {`Gaussian elimination'}, 1500 years before Gau\ss' birth and 1800 years after the middle-eastern seasonal flood prone countries had first used  the  Gaussian algorithm around 1800 BC.  
Gau\ss\  himself described the method as the  ``{common method of elimination}" in his papers and mathematicians then attached his name to it as an honor.
\vspace*{-1mm}

\subsection{\large The Genesis of Vectors and Matrices\ \ (3rd Epoch)}

To advance matrix computations further there was a need to conceptualize {coordinates} and {vectors} in space.\\
 In the 14th century AD, Nicole Oresme developed a system of {orthogonal coordinates} for describing Euclidean {space.}
  This idea was taken up by  Ren\'e Descartes in the 17th century and is familiar to all of us now under the concept of {Cartesian Coordinates}.
Thereby the world became ready for matrices and matrix computations in their own right.\\[1mm]
 In 1683 Gottfried Leibnitz in Germany and Seki Kowa in Japan both and unbeknownst to each other re-invented the {concept of a matrix} as a {rectangular array of coefficients} for studying linear equations.
  Leibnitz also used and suggested row elimination to simplify and find their solutions. These efforts enabled Gau\ss\ to repeat what  the Egyptians had done four millennia earlier:  He was asked to survey the lands of his ruler, the Archduke  George Augustus of Hanover and measure the size of this kingdom  inside Germany in the early 1800s. 
 Beginning in 1821 and 1823  Gau\ss\ , as Professor of Geodesy (and not of Mathematics) in G\"ottingen, would measure the angles and distances between many of the highest points there, such as the Brocken, the Inselsberg, 104 km apart,  and the hills around G\"ottingen and later expanded the surveys all the way to the North Sea. He did this multiple times, preferably  when the weather was clear. Thereby he set up  systems of linear equations with generally more equations than unknown due to repeated measurements on different days.\\[1mm]
   To solve these overdetermined and naturally `unsolvable'  system $Ax = b$, Gau\ss\ devised the {normal equation} $A^TAx = A^Tb$  and  solved it approximately.
 But the normal equations method eventually turned out to be  {numerically unsound}. It took over a century to find out why, the reason being that condition numbers multiply, see Olga Taussky (1950) \cite{Tauss}.
 
\subsection{\large Eigenvalues and the Characteristic Polynomial \ \  (4th Epoch)}

As  {differential operators} and {matrices} were beginning to be investigated and dealt with  by the early 1800s, their connections and similarities were slowly recognized in the mathematics world. \\[1mm]
The replication of certain functions $f \neq 0$ by a given differential operator $\cal A $ was noticed first and  became the subject of studies. What were the {functions $f$} for which {${\cal A} f = \alpha f$} for some scalar $\alpha$? How could they be found from $\cal A$, what about $\alpha$?\\[1mm]
 In 1829 Augustin Cauchy \cite[p. 175]{C1829} began to view the erstwhile {`eigenvalue equation'} ${\cal A} f = \alpha f$} as a {`null space equation'}, namely {${\cal A} f - \alpha f = 0$} or $({\cal A} - \alpha \ id)  f= 0$ for the identity operator with $id\ f = f$ for all $f$.
Complete knowledge of the {eigenvalues} $\alpha$ and {eigenfunctions} $f$ of a differential  operator $\cal A $ allowed for a simple sum representations of  the general solution of the  linear differential  equation described by $\cal A$. Thus Cauchy's `null space equation'  became essential for determining the behavior of systems governed by linear differential equations.\\[1mm]
    Cauchy's knowledge of and interest in determinants (think Cauchy-Binet Theorem) then led him to define the {`characteristic polynomial'} of a square matrix $A$  in 1839 \cite[p. 827]{C1839} as {$f_A(\alpha) = \det(A-\alpha I)$} and thereby he initiated renewed studies in polynomial root finding  algorithms in the hope of obtaining analogous diagonalization results for linear matrix times vector products. And the search for polynomial rootfinders was on. By modern day hindsight, reducing the eigenvalue problem from an $n^2$ data problem of the entries of a matrix $A_{n,n}$ to one of the $n+1$ coefficients of its characteristic polynomial is data compression and therefore it was doomed to fail. But  that remained unrealized by  the mathematics community for more than 100 years.\\[1mm]
     James Sylvester  finally gave the tableau concept of matrices its name {'matrix'} in 1848 or 1850. And after roughly 2 decades  \emph{1829 --$>$ 1839 --$>$ 1848/50}, 
  the first century  of {Matrix Theory} or theoretical {Linear Algebra} had begun.\\[-3mm]
 
  Back to { matrix computations}:\\
Cauchy's idea  led mathematicians  to try and compute the characteristic polynomials of matrices and find their {roots}  in order to understand the {eigen-behavior} of matrices.
 We still teach many concepts and lessons today that are based on the 'characteristic polynomial' $f_A(x)$ of a matrix $A$. Why, we should ask ourselves. 
Because unfortunately studying 'characteristic polynomials' in place of matrices has turned out  be a costly {dead end} for computational and applied mathematics: 
 In the century and a half that followed Cauchy's work, more than 4,000 papers on {computing the roots of polynomials} were published, together with 200 to 300 books on the subject, bringing us many algorithms, all of which failed more often than not.
  Many illustrious careers and schools of mathematics were founded based on this  unfortunate and  ever elusive goal.\\[1mm]
 During the same period, 2-D hand crank computing machines were invented and built to effect long number multiplications and divisions. First by Charles Babbage, then as commercial geared adding machines that were still being used in office work well into the 1960s.
      These worked as {two-dimensional abaci} of sorts.
 But eventually digital (at first punch-card fed) computers became the tools of our computational trade in science, in engineering, in business, in GPS, in Google, in social media,  large data, automation, etc, etc.
 
 But how could we or would we find {matrix eigenvalues} accurately?  A turn-around, a new method, a new computational Epoch was needed.\ \  
From where, by whom, and how ...?

\vspace*{-2mm}
\subsection{\large Iterative Matrix Algorithms \ \ (5th Epoch)}
To move us forward, it appears that {matrix methods} themselves might have to be developed that would solve the matrix-intrinsic  eigenvalue problem by themselves. 
But before that was possible there were further  unfortunate `detours'. \\[1mm]
Carl Friedrich Gau\ss\ - in his doctoral thesis of 1799 \cite{G1799} - had disproved all earlier attempts to establish the {Fundamental Theorem of Algebra}, i.e.,  that all polynomials over the reals  numbers can be factored into as many real or complex conjugate factors as their degree says. His thesis then included the first complete and correct proof of the Fundamental Theorem of Algebra.\\[1mm]
 In 1824 Niels Abel \cite{A1824} showed that the roots of some {5th degree polynomials} can not be found by using radical expressions of their coefficients; Gauss never opened or read the submitted paper and thus in fact rejected it knowingly on the `grounds' that God would not have complicated the World thus ... for us. 
 Abel published his result privately, a broken man. 
 \'Evariste Galois \cite{G1830a,G1830b} extended Abel's result in 1830 by giving group theoretic conditions for polynomials to be {solvable by radicals}; the extended paper (introducing Galois Theory) was originally rejected and appeared only posthumously in 1846 \cite{G1846}.\\[1mm] 
 Cautioned by these `rejected' inconvenient results, the  polynomial approach  to matrix eigenvalue computations  could have  been shunned by clearer heads early on; but the {`dead end' determinants and characteristic polynomial roots road} was taken instead for more than a century. Note that Cauchy's matrix result and most other fundamental matrix results from the  19th century were formulated in terms of determinants and only in the mid 20th century did the term 'matrix' appear in matrix theoretical article and book titles. \\[1mm]
A {matrix based approach} to the eigenvalue problem nowadays starts from the simple fact that for any $n$ by $n$ matrix $A$ and any $n$ vector $b$, the sequence of vector iterates $ b, Ab, A^2b, ..., A^kb,..., A^nb$ contains $n+1$ vectors in $n$-space which makes these vectors linearly dependent. Their linear dependency then leads to an $n$th degree polynomial $p_b(A)$ that sends $b$ to zero. The vanishing polynomial for any $b$ turns out to always be a factor of the characteristic polynomial of $A$ and  it can be found by Gaussian elimination rather than using determinants. \\[1mm]
 The same idea shows that {vector iteration} converges for every starting vector $b \neq 0$ and any given square matrix $A$ and this has lead engineers in the early 20th century to construct  iterative matrix algorithms that could solve linear equations and the matrix eigenvalue problem. 
Iterative matrix algorithms actually do go back further, such as to  the Jacobi (1839) \cite{Jaco}, Gau\ss\ -Seidel (1874) , and various SOR  methods that are designed to solve linear systems iteratively. These use matrix splittings of $A$ rather than vector iteration.\\[1mm]
 Alexei Krylov (1931) \cite{Kryl} introduced and studied the {vector iteration subspaces\ span$\{b, Ab, ..., A^kb\}$} in their own right. 
  Following his ideas, large sparse matrix systems are nowadays treated iteratively in so called Krylov-based methods, both to solve linear equations and to find matrix eigenvalues. 
  Standard widely used {Krylov type iterative matrix algorithms}  carry the names of Steepest Descent, Conjugate Gradient, Arnoldi (1951) \cite{Arno}, Lanczos (1950) \cite{Lancz}. Others are called  GMRES, BICGSTABLE, QMR, ADI etc. Most Krylov type methods are matrix and problem specific and they are mostly used for huge sparse and structured matrices where direct or semi-direct methods cannot be employed due to their high computational and storage costs. Krylov methods generally rely on preconditioners $M$ for linear systems $Ax = b$ that shift the spectrum of $M^{-1}A$ for faster convergence and they thrive on incomplete matrix splittings etc. Typically they give only partial results. Who needs to know all the  million eigenvalues of a million by million matrix model. Krylov methods can be tuned to give  information where it is needed  for the underlying system.\\[-6mm]
  
  \subsection{\large Francis Algorithm and Matrix Eigenvalues \ \  (6th Epoch)}

The Second World War and post Second World War periods were filled with innovations: \\
 The atomic era had begun, as well as rocket science; commercial air flight became popular; digital computers were being developed, first as valve machines and later transistorized. Supersonic speeds were realized, Computer Science was developed, etc.
But there were many crashes and disasters with the new technologies: Commercial aircraft (Super-Constellation, Convair, ...) and military ones (Starfighter ...) would crash weekly around the globe; and newly built suspension bridges would collapse in strong winds. \\[1mm]
The crux of the matter was that while {matrix models} of the underlying mechanical systems could readily be made using the laws of physics and mechanics, no one could reliably compute their eigenvalues. 
Engineers could not test their designs for {eigen-modes in the right half plane}! And Krylov methods were unfortunately not sufficient  for testing for eigenvalues in a half plane.\\[1mm]
 If a {matrix model of a mechanical or electrical or ... structure}   has right-half plane eigenvalues $\lambda$ then - upon proper excitation - there would be a eigen component of the ever increasing form $e^{\lambda t} \to \infty$ as $t \to \infty$ that resonates and self-amplifies  inside the structure itself. This  then leads to ever increasing {destructive vibrations} and ultimate failure.
 The aircraft 'flutter problem' was discovered during World War 2. In England during WW2,  Gershgorin circles that contain all of a  system's eigenvalues were drawn  out in the complex plane by rather primitive valve computers and checked to ascertain  system stability. \\[1mm]
    The general {matrix eigenvalue problem} was finally solved independently  and similarly by {John Francis} in London and by {Vera Kublanovskaya} in Russia nearly simultaneously around 1960. 
Francis' (or the QR) algorithm \cite{Franc1, Franc2} is based on Alston Householder's idea to try and solve matrix problems by {matrix factorizations}. Francis' method is an orthogonal  subspace projection method and it works  differently  than the Krylov based methods which solve a given matrix  eigenvalue problem by projecting  onto a Krylov subspace that is derived from and suitable for  $A$.\\[1mm]
 A 'divide and conquer' matrix factorization strategy  was first employed by Heinz Rutishauser (1955, 1958) \cite{Ruti1, Ruti2} in his {LR matrix eigenvalue algorithm}: If one can {factor} $A = LR$ into the product of  a lower and an upper triangular  matrix $L$ and $R$ as $A = LR$ and if $L$ is invertible, then for the reverse order product $A_1 = RL$ we have $A_1 = L^{-1} AL$ since $R = L^{-1}A$. 
If $A_1$ again allows an LR factorization   $A_1 = L_1 R_1 $ with $L_1$  nonsingular, then by {reverse order multiplying} we obtain 
  $$A_2 = L_1^{-1} A_1 L_1 = L_1^{-1} L^{-1} AL L_1$$
   and so for the sequence of likewise constructed matrices $A_i$ for $i=3,...$ if the respective LR factorizations are possible at each stage $i$.   In this case the iterates $A_i$  clearly remain similar to the original  matrix $A$ and thus the iterates all have the same eigenvalues as $A$. Note, however, that if for example the (1,1) entry  $A(1,1)$ is zero in $A$, then there exists no un-pivoted LR factorization for A and the method breaks down since pivoting is a one sided matrix process and not a similarity. Therefore Rutishauser's method is only applicable to a limited set of matrices $A$ for which every LR iterate $A_i$ allows an un-pivoted LR factorization. 
  Rutishauser had noted that if  LR factorizations are possible for all iterates $A_i$, then the $A_i$  become nearly {upper triangular} for large $i$ with $A$'s {eigenvalues on the diagonal}. (Very loosely said.)\\[1mm]
 John Francis was very interested in the flutter problem at the time when, by chance, someone dropped Rutishauser's 1958 LR paper \cite{Ruti2} on his desk at the CRDC in London. \ \emph{(In my  interview with John Francis in 2009 \cite{GolUhinterv}, he did not know who that might have been.)} \ 
 Francis was aware through contacts with Jim Wilkinson of the {backward stability} of algorithms that involve {orthogonal matrices} Q. So rather than using Gaussian elimination matrices $L$, Francis experimented with orthogonal $ A = QR$ factorizations. 
 At roughly the same time Vera Kublanovskaya worked on an LQ factorization of $A$ as $A = LQ$ and subsequent reverse order multiplications \cite{Kubla} in Leningrad, Russia, that would also compute the eigenvalues of $A$.\\[1mm]
 Rutishauser had  observed convergence speed-up for his LR method when replacing $A$ by $A - \alpha I$, i.e., shifting. 
Hence   Francis  experimented  with shifts for QR and then established the {'Implicit Q Theorem'} \cite{Franc2} in order to circumvent computing {eigenvalues of real matrices} over the complex numbers. 
Implicit shifts  also avoid  {rounding errors} that would be introduced by explicit shifts. 
   Francis' second paper (1962) \cite{Franc2} also contains a fully computed flutter matrix problem of size 24 by 24. The eigenvalues of such 'large' problems had never before been computed successfully.\\[1mm]
Francis' Implicit Q Theorem then allowed Gene Golub and Velvel Kahan (1965) \cite{GolKah} to compute {singular values} of large matrices for the first time and this application later spawned the original Google search engine and brought us - in a way - into the {internet age}.\\[1mm]
 In 2002 the {Multishift QR Algorithm} was developed by Karen Braman, Ralph Byers, and Roy Mathias \cite{BraBM1, BraBM2}. It relies on {subspace iteration} and extends  Francis' QR and combes it with Krylov like methods.
 This extension allows us today to compute the complete eigen and singular value structure of dense matrices of sizes up to 10,000 by 10,000 economically on laptops.\\[-3mm]

What is being missed today computationally? What Epoch(s) might come next? Why and how?

\subsection{{\large Two New Epochs Ahead \ \  (7th and 8th Epochs)} \ -- yet to come}

Two new Epoch generatinging impulses have become visible on the matrix computational horizon of today:\\[2mm]
{\bf (A) } One expands our computational abilities from static problem solving algorithms to  {real time methods for time-varying problems}, and\\[2mm]
 {\bf (B) } the other 
involves {computer hardware} advances.\\

{\bf 1.7.1 \ Time-varying Problems and Real-time Solvers  \ \emph{( Epoch 7)} :}\\[4mm]
Our current best numerical codes can solve static problem very well; that is what they are designed for. \\[1mm]
As we begin to rely more and more on time-dependent sensoring and on robotic manufacture, we need to learn how to solve our erstwhile {static equations}, but now  in real time and {with time-varying coefficients}, and preferably accurately as well. 
It seems quite alluring to try and solve a time-varying problem by using the static time-dependent inputs at each instance statically. 
 But such a naive solution cannot suffice since at the next time step, whose `solution' we have just computed `statically', the problem parameters have already changed and thus our `static' solution solves a completely different problem, which -- unfortunately -- has little value. \hfill \emph{If any at all}.\\
 
{\bf 1.7.2 \ Computer Hardware \ \emph{( Epoch 8 )} :}\\[2mm]
Since the earliest electronic computing devices of the 1940s, all our  computers have worked as giant and embellished Turing machines with logic gates, switches and memory that rely only two numerical states: 0 and 1 or on or off. Hence  all our computer data is stored and manipulated as sequences of {0 and 1}.\\[1mm]
Lately our computing ability has  come up against the limits of storing and working with data and processors that  can only deal with zeros and ones. 
 Our {processing speeds}   have not advanced significantly over the last couple of years; we are still stuck with 3 - 4 GHz processors. To alleviate this bottleneck, chip makers have created multi-processor chips and software firms have introduced better and quicker software and operating systems, but the basic processor speed has not budged much.\\[1mm]
At this time computer scientist and manufacturers are  trying to overcome this 0 - 1 bottle-neck by replacing our 0 - 1 {processors, chips, memories and transistors} by improved transistors and chips that can store and process {multi-states}, such as 0-1-2-3-4 or 0-1-2-3-4-5-6-7-8 or even higher numbered data representations. 
 This could lead us to another computing  sea change bringing us into a new computational Epoch via hardware. 
And further out on the horizon lies the possibility of having {infinitely many quantum states based computers}.

\section{\Large On Neural Network Methods  \ \ { \emph{(Epoch 7}} \ \ \ \normalsize{already under way}\Large)}

The last century brought us valuable tools to solve most  static problems that involve matrices.\\[2mm]
 Our current numerical matrix tools can solve  static  linear equations,  matrix equations such as Sylvester or Lyapunov equations, eigenvalue problems and generalizations of these, both of the  dense or structured and of the solvable or unsolvable kind; likewise we can solve static optimization problems of all sizes and for nearly all structured matrices  and for most static applications. \\[1mm]
But what do do with such problems when their coefficients are time-varying or time dependent?\\
In numerical computations, there has always been a  see-saw between  models that resulted in derivative inspired differential equations and in linear algebra based matrix equations. Their respective computational advantages differed from problem to problem. Neural networks (NN) are an amalgam of matrix methods and differential methods using a mixture of both. NN methods are designed to solve time-varying dynamical systems. Numerical methods for time-varying dynamical systems first came about in the 1950s and subsequently have gained strength in the 1980s and 1990s and beyond, see the introduction in Getz and Marsden \cite{GM97} for example. To solve a time varying equation $f(X(t)) = G(t)$, these studies start  from the error equation $E(t) = f(X(t)) - G(t)$ and stipulate exponential decay of the error function $E(t)$ by trying to solve 
\begin{equation}\label{erroreq}
 \dot E(t) = - \lambda E(t)
 \end{equation}
for a positive decay constant $\lambda$. There are essentially three ways to go about solving the error function differential equation (\ref{erroreq}): homotopy methods, gradient methods and ZNN Neural Network methods.

\subsection{\large A Neural Network approach to solve {time-varying linear equations
 $\mathbf{A(t) x(t) = b(t)} $}}
Here $A(t)$ is a nonsingular time-varying $n$ by $n$ matrix and $b(t) \in \RR^n$ is a time-varying vector, respectively. Clearly the unknown solution $x(t)$ of the associated linear equation $A(t) x(t) = b(t)$ will  be time-dependent as well.\\[2mm] 
The first paper on Zhang Neural Networks (ZNN) was written by Yunong Zhang et al. in 2002, see \cite{Zh2002}. Today there are over 250 papers, mostly in engineering journals that  deal with time-varying applications of the ZNN method, either in hardware chip design for specialized computational tasks as part of a plant or   machine or for time-varying  simulation problems  in computer algorithms and codes. Unfortunately, the ZNN method and the ideas behind ZNN are hardly known today among  numerical analysis experts. The method itself  starts with using Suanjing's and Al Khwarizmi's  rule for reducing  equations that first appeared  1 1/2 millennia ago. This simple rule was also employed by  Cauchy in 1829 to transform the static matrix eigenvalue problem from $ A x = \lambda x $ to $ A x- \lambda x = 0$  and finally to  det$(A-\lambda I) = 0$.  For the time-varying linear equations problem Zhang's Neural Network method starts with $$ A(t) x(t) - b(t) = 0 \ .$$
ZNN then looks at the {error function} $$E(t) = A(t)x(t) - b(t) : \RR \to \RR^n$$
Note that standard static methods would look at an {error norm} $\|E\|$ for the error function.   Neural Networks do not, they study the {error function}  $E(t)$ instead. 
And they start with an implicit `{ideal wish}': 
 \emph{What could or  should we wish for $E$}? 
Time-varying computations would be near  ideal if their {error functions} were {decaying exponentially fast}.  
This is {impossible} to achieve (or even ask for) with our best static equations and problem solvers of the 21st century. For static numerical matrix methods backward stability is considered  best. \\
In all Zhang NN methods we stipulate that the {error function $E(t)$  decreases exponentially} fast over time to the zero function. 
This means that
$$ {\dot E(t) = -\lambda E(t)} \ \text{ for some chosen number } \lambda \gg 0, \  \text{  the decay constant}.$$
Note that for the time-varying linear equations problem we have 
$$\dot E(t) = \dot A(t)x(t)+ A(t)\dot x(t)-\dot b(t)\ . $$
This leads to the following {differential equation} for $x(t)$:
\begin{equation}\label{linequDE} 
A(t) \dot x(t) = -\dot A(t) x(t) + \dot b(t) -\lambda(A(t)x(t) - b(t))\ .
\end{equation}
And we have transformed the time-varying linear equations problem into an {initial  value differential equations problem} that needs to be solved for $t > 0$.  This is where the different dynamical systems methods split their ways. In Zhang Neural Networks, the continuous time differential equation  (\ref{linequDE}) is then discretized for $0 < t_j < t_{end}$ and the ensuing derivatives are approximated by high order difference quotients, with the one for the unknown $\dot x(t_j)$ being 1-step ahead and proven convergent, while the others such as for $\dot A(t)$ and $\dot b(t)$ can be backward difference formulas.  
How to proceed from (\ref{linequDE})  with solving $A(t)x(t) = b(t)$ via ZNN methods is still an open problem, especially  for large scale sparse or structured time-varying linear equations since the matrix $A(t_j)$ encumbers the unknown $\dot x(t_j)$ on the left hand side and there is no known 1-step ahead differentiation formula that can be used here.\\ 
The general {idea that underlies ZNN methods} for time-varying problems is to {replace repeated matrix computations} by solving {linear differential equations} and associated {initial value problems} for discrete instances $0 < t_j < t_{end}$ instead. \\[-5mm]

\subsection{\large A Zhang Neural Network approach to find  {time-varying generalized matrix inverses $\mathbf{ Y(t)}$ for time-varying full rank matrices $\mathbf{ B(t)}$ with 
 $\mathbf{ B(t)_{m,n} Y(t)_{n,m} = I_m} $}} 

This section is based on joint work with  Jian Li, Mingzhi Mao, and Yunong Zhang, see \cite{LMUZh2}.\\[1mm]

\noindent
{\bf{Continuous Problem Formulation} :}\\[-2mm]

\noindent
For an $m$ by $n$ real time varying matrix $B(t)$ of full rank $m$ with $m \leq n$ 
we form the matrix-valued error function
\begin{equation}\label{1}
E(t)=B(t)-Y^+(t)\in\mathbb{R}^{m\times n} 
\end{equation}
where the upper + sign always means 'generalized inverse'.
Then we use the  design formula 
\begin{equation}\label{2}
\dot{E}(t)=-\lambda E(t)
\end{equation}
with design parameter $\lambda>0$.
Based on \cite[Lemma 3]{LZh14} we have 
\begin{equation}\label{3}
\dot{Y}^+(t)=-Y^+(t)\dot{Y}(t)Y^+(t).
\end{equation}
And from (\ref{1}) and (\ref{3}) we obtain 
\begin{equation}\label{4}
\dot{E}(t)=\dot{B}(t)-\dot{Y}^+(t)=\dot{B}(t)+Y^+(t)\dot{Y}(t)Y^+(t).
\end{equation}
Combining (\ref{2}) and (\ref{4}), we then get 
\begin{equation}\label{5}
\dot{B}(t)+Y^+(t)\dot{Y}(t)Y^+(t)=-\lambda(B(t)-Y^+(t))
\end{equation}
And by right multiplying (\ref{5}) with $Y(t)$ we  have 
\begin{equation}\label{6}
(\dot{B}(t)+Y^+(t)\dot{Y}(t)Y^+(t))Y(t)=-\lambda(B(t)-Y^+(t))Y(t).
\end{equation}
With $m\leq n$, we have $Y_{m\times n}^+(t)Y_{n\times m}(t)=I_{m\times m}$ and  thus
\begin{equation}\label{7}
\dot{B}(t)Y(t)+Y^+(t)\dot{Y}(t)=-\lambda(B(t)Y(t)-I).
\end{equation}
The solution of a generalized-matrix inverse problem is not unique when $m< n$ and we only need to find a solution that satisfies (\ref{7}). Consequently  the continuous model can be represented as
\begin{equation}\label{8}
\begin{split}
\dot{Y}(t)=-\lambda (Y(t)B(t)Y(t)-Y(t))-Y(t)\dot{B}(t)Y(t)
\end{split}
\end{equation}
which agrees completely with \cite[formula (15), p. 317]{GM97}.
Substituting (\ref{8}) into (\ref{7}), we have 
\begin{equation}\label{9}
\dot{B}(t)Y(t)+Y^+(t)(-\lambda (Y(t)B(t)Y(t)-Y(t))-Y(t)\dot{B}(t)Y(t))=-\lambda(B(t)Y(t)-I),
\end{equation}
which we rewrite as
\begin{equation}\label{10}
\dot{B}(t)Y(t)+(-\lambda (Y^+(t)Y(t)B(t)Y(t)-Y^+(t)Y(t))-Y^+(t)Y(t)\dot{B}(t)Y(t))=-\lambda(B(t)Y(t)-I).
\end{equation}
With $Y_{m\times n}^+(t)Y_{n\times m}(t)=I_{m\times m}$, we have
\begin{equation}\label{11}
\dot{B}(t)Y(t)+(-\lambda (B(t)Y(t)-I)-\dot{B}(t)Y(t))=-\lambda(B(t)Y(t)-I).
\end{equation}
Thus model (\ref{8}) satisfies (\ref{7}) and its solution solves the time-varying generalized matrix inverse problem. \\[-2mm]

\noindent
{\bf Zhang Neural Network Discretization :}\\[-2mm]

Given a sequence of rectangular matrices $B_j$ at  time instances $t_j \leq t_k$ we want to find  the discrete time-varying generalized matrix inverse $Y_{k+1}$ of $B_{k+1}$ on each computational time interval $[k\tau,(k+1)\tau)\subseteq[0,t_{\textnormal{f}}]$ so that 
\begin{equation}\label{d.p.pseudoinverse}
B_{k+1}-Y^+_{k+1}={0}.
\end{equation}
Here $B_{k+1}=B(t_{k+1})=B((k+1)\tau)\in\mathbb{R}^{m\times n}$ is a time-varying full rank equidistant matrix sequence,  $m\leq n$ and $Y_{k+1}\in\mathbb{R}^{n\times m}$ is unknown.  $Y_{k+1}$ needs to be computed in real-time for each time interval $[k\tau,(k+1)\tau)\subseteq[0,t_{end}]$. Here   the matrix operator $^+$ denotes the generalized inverse of a matrix and ${0}\in\mathbb{R}^{m\times n}$ is the zero matrix. Besides, $k=0,1,\cdots$ denotes the updating index, $t_{end}$ denotes the task duration and $\tau$ denotes the constant sampling gap of the  time-varying matrix sequence $B_{k+1}$. For $m>n$, the procedure  is similar.

Note that we must obtain each $Y_{k+1}$ at or before time $t_{k+1}$ for real-time calculations while the actual value of $B_{k+1}$ is unknown before $t_{k+1}$. Thus we cannot obtain the solution  by calculating $Y_{k+1}=B^+_{k+1}$. To obtain $Y_{k+1}$ in real-time, we must develop a model based on the available information from before $t_{k+1}$ such as that in $B_j$, $Y_j$ and $Y_{j-1}$ for $j \leq k$  instead of  unknown information such as  $B_{k+1}$.

To obtain a discrete-time model that solves the original discrete time-varying generalized matrix inverse problem (\ref{d.p.pseudoinverse}), we need to discretize the continuous model (\ref{8}). First we use the  conventional 1-step forward Euler  formula 
\begin{equation}\label{e.f.formula}
\dot x(t_{k})\approx\frac{x(t_{k+1})-x(t_{k})}{\tau}
\end{equation}
with  truncation error of order $O(\tau)$.
Based on  formula (\ref{e.f.formula}) we approximate
\begin{equation}\label{e.f.formula.A}
\dot Y(t_{k})=\frac{Y(t_{k+1})-Y(t_{k})}{\tau} \ 
\end{equation}
and use this equation  to discretize the continuous  model (\ref{8}) as follows
\begin{equation}\label{e.d.m.pseudoinverse.k}
\begin{split}
Y_{k+1}=-h (Y_kB_kY_k-Y_k)-\tau Y_k\dot{B}_kY_k+Y_k .
\end{split}
\end{equation}
Here $h=\tau\lambda$. In most real-world applications, information of the first-order time derivatives, i.e., the value of $\dot{B}_k$ may not be explicitly known for the discrete time-varying generalized matrix inverse problem (\ref{d.p.pseudoinverse}). If this is so, the value of $\dot{B}_k$ can be approximated by a backward finite difference formula. To assure the accuracy and simplicity of the discretized  model, the truncation error of the backward finite difference formula for $\dot{B}_k$ should be near equal to that of the 1-step-ahead finite difference formula that approximates $\dot{Y}_k$. Thus $\dot{B}_k$ in (\ref{e.d.m.pseudoinverse.k}) should best be approximated by Euler's backward finite difference formula 
\begin{equation}\label{e.b.f}
\dot b_k\approx \frac{b_k-b_{k-1}}{\tau},
\end{equation}
because the truncation error order $O(\tau)$ of formula (\ref{e.b.f}) equals that of  formula (\ref{e.f.formula}). Thus  we approximately have
\begin{equation}\label{e.b.f.A}
\dot B_k= \frac{B_k-B_{k-1}}{\tau}.
\end{equation}
Then we combine equation (\ref{e.d.m.pseudoinverse.k}) with equation (\ref{e.b.f.A}) and the Euler discrete model becomes
\begin{equation}\label{e.d.m.pseudoinverse}
\begin{split}
Y_{k+1}=-h (Y_kB_kY_k-Y_k)- Y_k(B_k-B_{k-1})Y_k+Y_k.
\end{split}
\end{equation}
Note that the truncation error of the  discrete model (\ref{e.d.m.pseudoinverse}) is of order $\mathbf{O}(\tau^2)$ where the symbol $\mathbf{O}(\tau^2)$ denotes a matrix in which each entry is of order ${O}(\tau^2)$. This  model  uses only  present or past information  of $B_k$, $B_{k-1}$ and $Y_k$ and solves for $Y_{k+1}$. Thus  $Y_{k+1}$ can be calculated during the time interval $[t_{k},t_{k+1})$ and if  $Y_{k+1}$ can be computed quickly enough  in real-time it will be ready when time instant $t_{k+1}$ arrives. 

Higher accuracy 1-step ahead formulas exist for discrete models, namely
\begin{equation}\label{4.1.formula}
\dot x(t_{k})\approx\frac{2x(t_{k+1})-3x(t_{k})+2x(t_{k-1})-x(t_{k-2})}{2\tau}
\end{equation}
and
\begin{equation}\label{4.2.formula}
\dot x(t_{k})\approx\frac{6x(t_{k+1})-3x(t_{k})-2x(t_{k-1})-x(t_{k-2})}{10\tau} .
\end{equation}
Both have truncation errors of  order $O(\tau^2)$. For simplicity we only consider formula (\ref{4.1.formula})  and call it the  4-IFD formula  because 4 instants in time are used to approximate the first-order derivative of $x(t_{k})$. When we employ the 4-IFD formula (\ref{4.1.formula}) inside the our continuous  model (\ref{8}) we obtain
\begin{equation}\label{41.d.m.pseudoinverse.k}
\begin{split}
Y_{k+1}=-h (Y_kB_kY_k-Y_k)-\tau Y_k\dot{B}_kY_k+\frac{3}{2}Y_k-Y_{k-1}+\frac{1}{2}Y_{k-2}.
\end{split}
\end{equation}
Next we use the three-instant backward finite difference formula
\begin{equation}\label{3.b.f}
\dot b_k\approx\frac{3b_k-4b_{k-1}+b_{k-2}}{2\tau}
\end{equation}
with error order $O(\tau^2)$ to approximate the value of $\dot{B}_k$ in equation (\ref{41.d.m.pseudoinverse.k}). Then the 4-IFD-type discretized  model becomes
\begin{equation}\label{41.d.m.pseudoinverse}
\begin{split}
Y_{k+1}=-h (Y_kB_kY_k-Y_k)- Y_k\left(\frac{3}{2}B_k-2B_{k-1}+\frac{1}{2}B_{k-2}\right)Y_k+\frac{3}{2}Y_k-Y_{k-1}+\frac{1}{2}Y_{k-2}.
\end{split}
\end{equation}
Its truncation error is of order $\mathbf{O}(\tau^3)$. Similar to the Euler based discrete model (\ref{e.d.m.pseudoinverse}), the  4-IFD-type discrete model  uses only  present and past information such as  $B_k$, $B_{k-1}$, $B_{k-2}$, $Y_k$, $Y_{k-1}$ and $Y_{k-2}$ to solve for $Y_{k+1}$. Thus it  also satisfies the requirements for real-time computation.\\[-2mm]

\noindent
{\bf{A 5 instants finite difference  formula} :}\\[-2mm]   

Any  usable finite difference formula for discretizing the continuous  model (\ref{8}) must satisfy several restrictions. It must be 1-step ahead for $\dot x$, i.e., approximate $\dot{x}(t_k)$ by using $x_{k+1}, \ x_{k}, \ x_{k-1}$ and earlier data, and it must be $0$-stable and convergent. However,  1-step ahead finite difference formulas do not necessarily generate  stable and convergent discrete models, see e.g.,  \cite{CAM2014JL,CAM2013CY}.

Here is a new 1-step ahead finite difference formula with  higher accuracy than the Euler  and 4-IFD formulas. It will be used to generate a stable and convergent discrete model that  finds time-varying generalized matrix inverses more accurately in real-time.\\[2mm]
{\bf Theorem 1} \
\emph{ The 5-IFD formula
\begin{equation}\label{new.formula}
\dot{x}(t_k)\approx\frac{8x(t_{k+1})+x(t_k)-6x(t_{k-1})-5x(t_{k-2})+2x(t_{k-3})}{18\tau}
\end{equation}
has  truncation error order $O(\tau^3)$.}\\[2mm]
\noindent
The proof relies on four Taylor expansions that use $x(t_{k+1})$ and $x(t_{k-1})$ through $x(t_{k-3})$ around $x(t_k)$ and clever linear combinations thereof.\\[1mm]
The new 1-step ahead discretization formula (\ref{new.formula}) then leads to the 5 instants discrete model 
\begin{equation}\label{5.d.m.pseudoinverse}
\begin{split}
Y_{k+1}=&-\frac{9}{4}h (Y_kB_kY_k-Y_k)-\frac{9}{4} Y_k\left(\frac{11}{6}B_k-3B_{k-1}+\frac{3}{2}B_{k-2}-\frac{1}{3}B_{k-3}\right)Y_k\\
&\qquad -\frac{1}{8}Y_k+\frac{3}{4}Y_{k-1}+\frac{5}{8}Y_{k-2}-\frac{1}{4}Y_{k-3}.
\end{split}
\end{equation}
 which has a  truncation error  of order $\mathbf{O}(\tau^4)$.\\[2mm]
{\bf Theorem 2}\label{p.0stable} \
\emph{The 5 instants discrete model (\ref{5.d.m.pseudoinverse}) is 0-stable.}\\[2mm]
\noindent
The multistep formula of the 5 instants discrete model  time-varying generalized matrix inverses has the  characteristic polynomial
\begin{equation}\label{characteristic.polynomial}
P_4(\theta)=\theta^4+\frac{1}{8}\theta^3-\frac{3}{4}\theta^2-\frac{5}{8}\theta+\frac{1}{4}
\end{equation}
with four distinct roots $\theta_{1,2}=-0.7160\pm 0.5495\textnormal{i}$, $\theta_3=0.3069$ and $\theta_4=1$ inside the complex unit circle, making this model 0-stable.\\[-2mm]

\noindent
{\bf Numerical Examples :}\\[-2mm]

{\it Example 1.} Consider the  discrete time-varying generalized matrix inverse problem
\begin{equation}\label{example1}
B_{k+1}-Y^+_{k+1}={0},~\textnormal{with} ~B_{k}=
\begin{bmatrix}
\sin(0.5t_{k})&\cos(0.1t_{k})&-\sin(0.1t_{k})\\
-\cos(0.1t_{k})&\sin(0.1t_{k})&\cos(0.1t_{k})
\end{bmatrix}.
\end{equation}
\newpage

\begin{figure}\vspace*{-14mm}\centering
\psfrag{k}[c][c][1.6]{$k$}%
\psfrag{5ifd}[c][c][0.5]{~~~~~~~~~5 instants} 
\psfrag{4ifd}[c][c][0.5]{~~~~~~~~~4 instants} 
\psfrag{euler}[c][c][0.5]{~~~~~Euler formulas} 
\psfrag{eee}[c][c][0.6]{~~~~~~~~$\|B_{k+1} Y_{k+1}-I\|_{\textnormal{F}}$}%
\vspace*{-6mm}
\subfigure[With $\tau=0.1$ s]{\includegraphics[width=0.325\columnwidth]{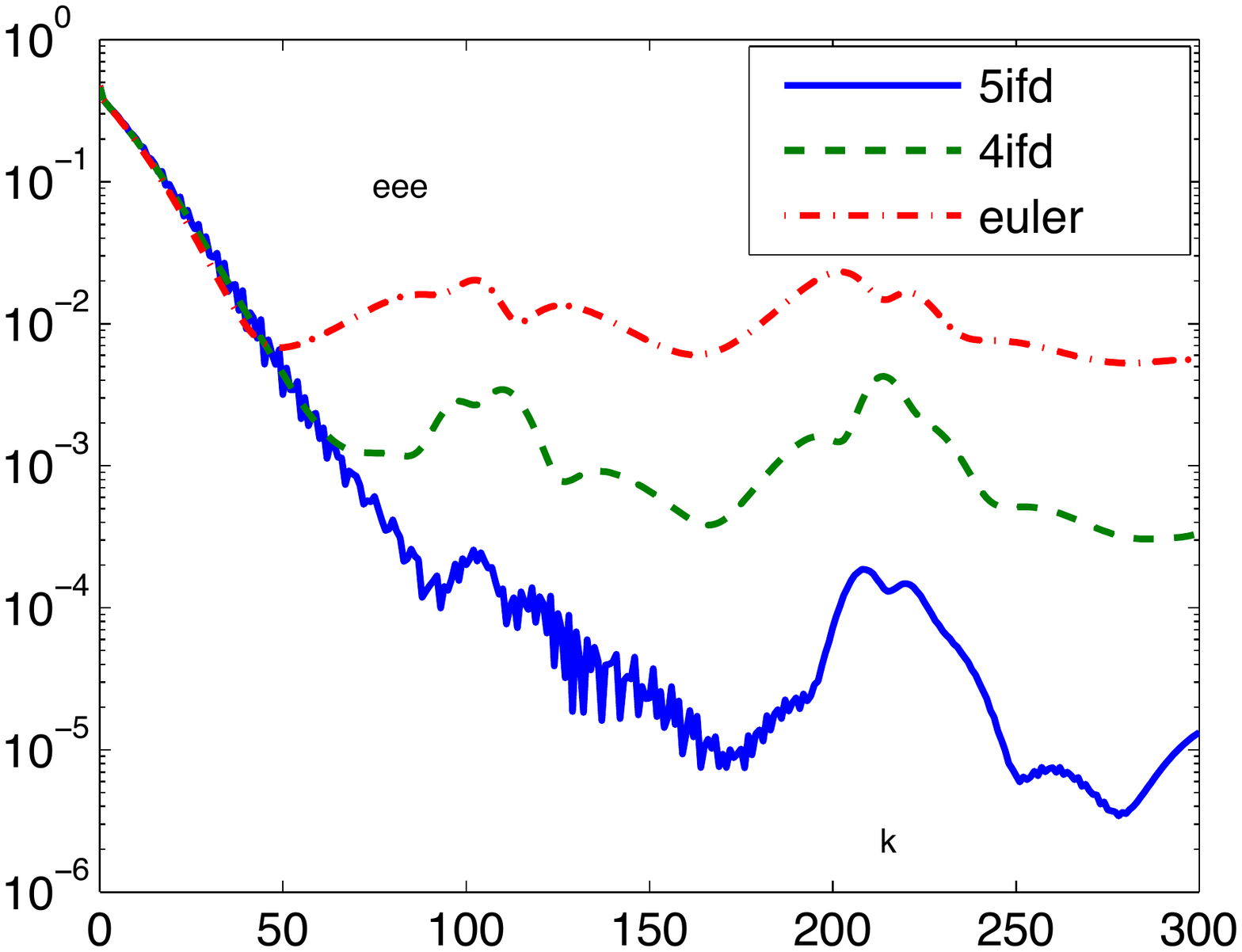}}
\subfigure[With $\tau=0.01$ s]{\includegraphics[width=0.325\columnwidth]{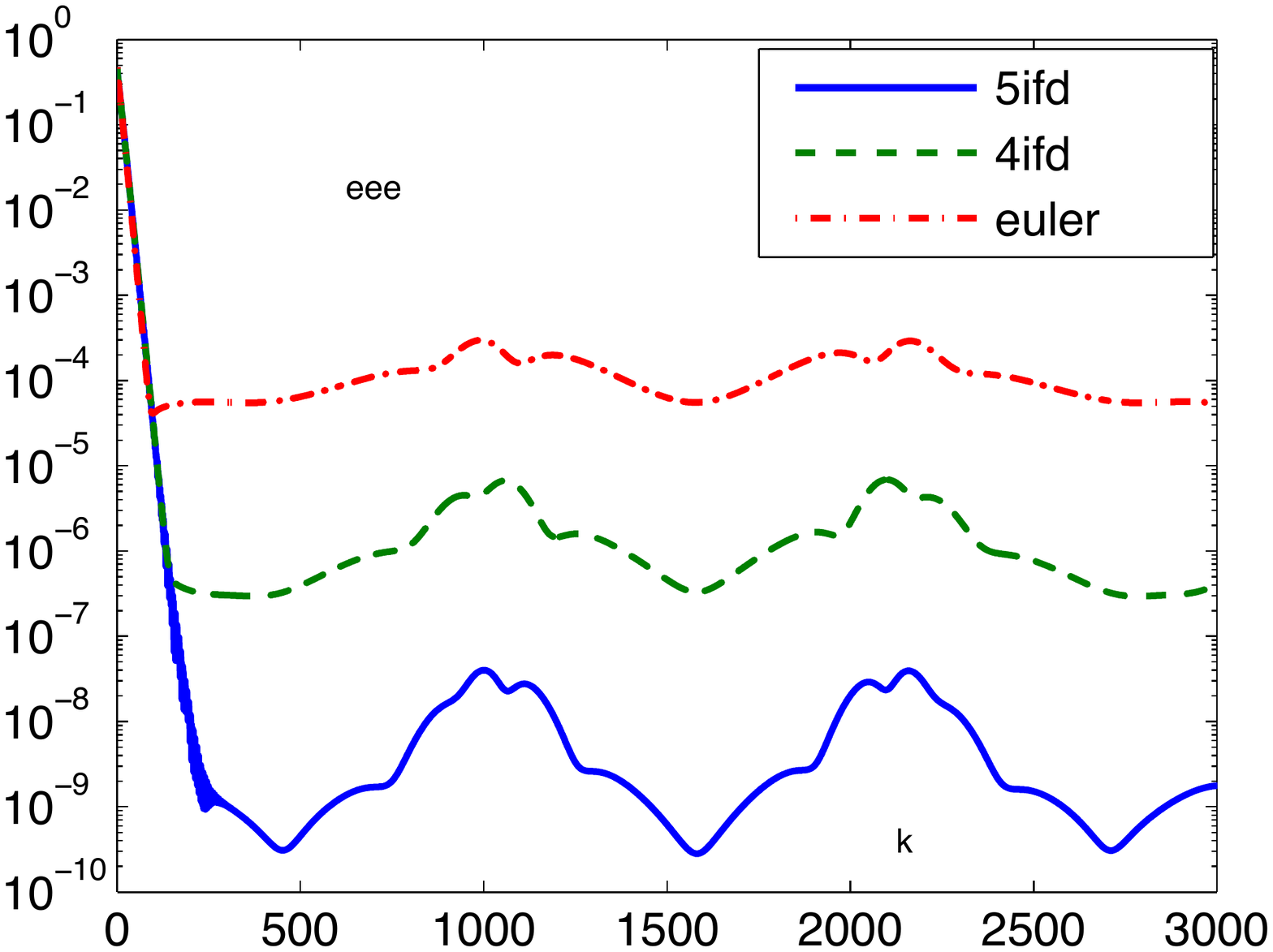}}
\subfigure[With $\tau=0.001$ s]{\includegraphics[width=0.325\columnwidth]{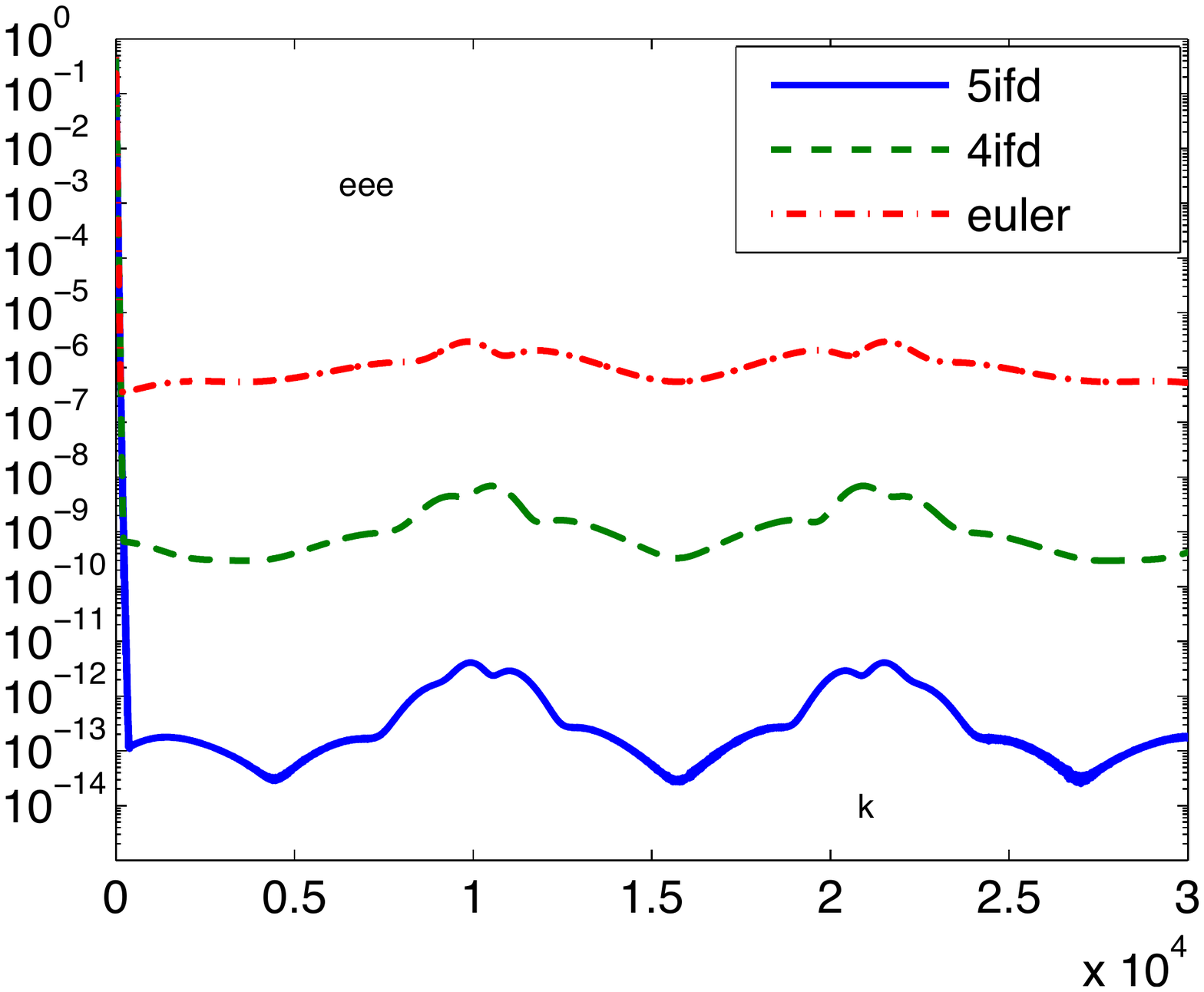}}
\caption{Typical residual errors generated by the 5 instants , the 4 instants   and the Euler 
formulas  with different  sampling gaps $\tau$ when solving the discrete time-varying generalized-matrix-inverse problem (\ref{example1}) for $t_{end}= 30$ s, and $h = 0.1$.} 
\label{errors.p}
\end{figure}

\vspace*{0mm}

{\it Example 2.} Here we consider the discrete time-varying matrix inversion problem
\begin{equation}\label{example2}
A_{k+1}X_{k+1}=I ~~~ \textnormal{with} ~A_{k}=
\begin{bmatrix}
\sin(0.5t_{k})+2&\cos(0.5t_{k}) \\
\cos(0.5t_{k})&\sin(0.5t_{k})+2 \\
\end{bmatrix}.
\end{equation}
\vspace*{-8mm}
\begin{figure}\vspace*{-10mm}\centering
\psfrag{k}[c][c][0.7]{$k$}%
\psfrag{x11}[c][c][0.7]{$X_{1,1}$}%
\psfrag{x12}[c][c][0.7]{$X_{1,2}$}%
\psfrag{x21}[c][c][0.7]{$X_{2,1}$}%
\psfrag{x22}[c][c][0.7]{$X_{2,2}$}%
\subfigure[Profiles of the (1,1) entry  $X_{1,1}$]{\includegraphics[width=0.4\columnwidth]{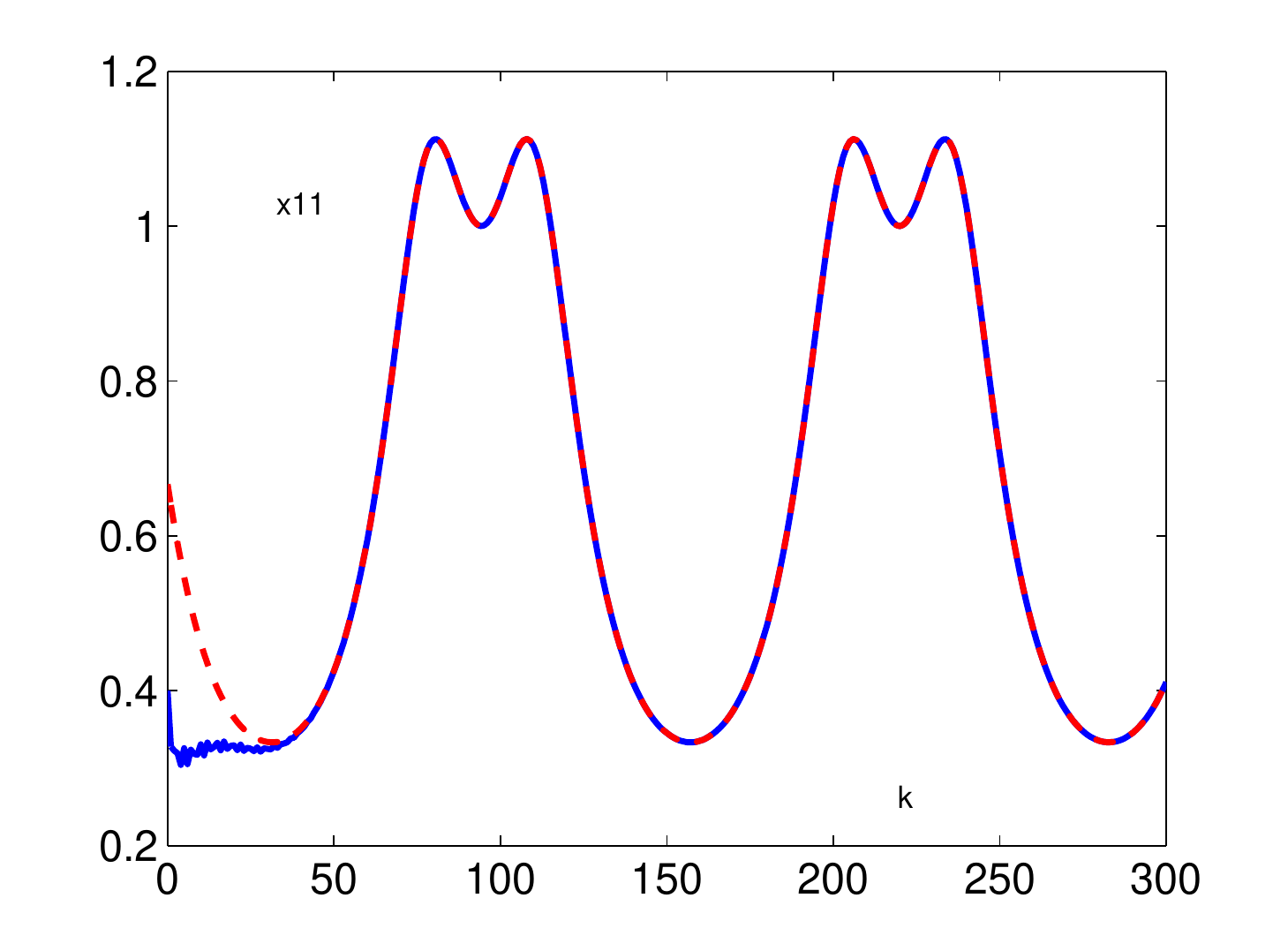}}
\subfigure[Profiles of the (1,2) entry $X_{1,2}$]{\includegraphics[width=0.4\columnwidth]{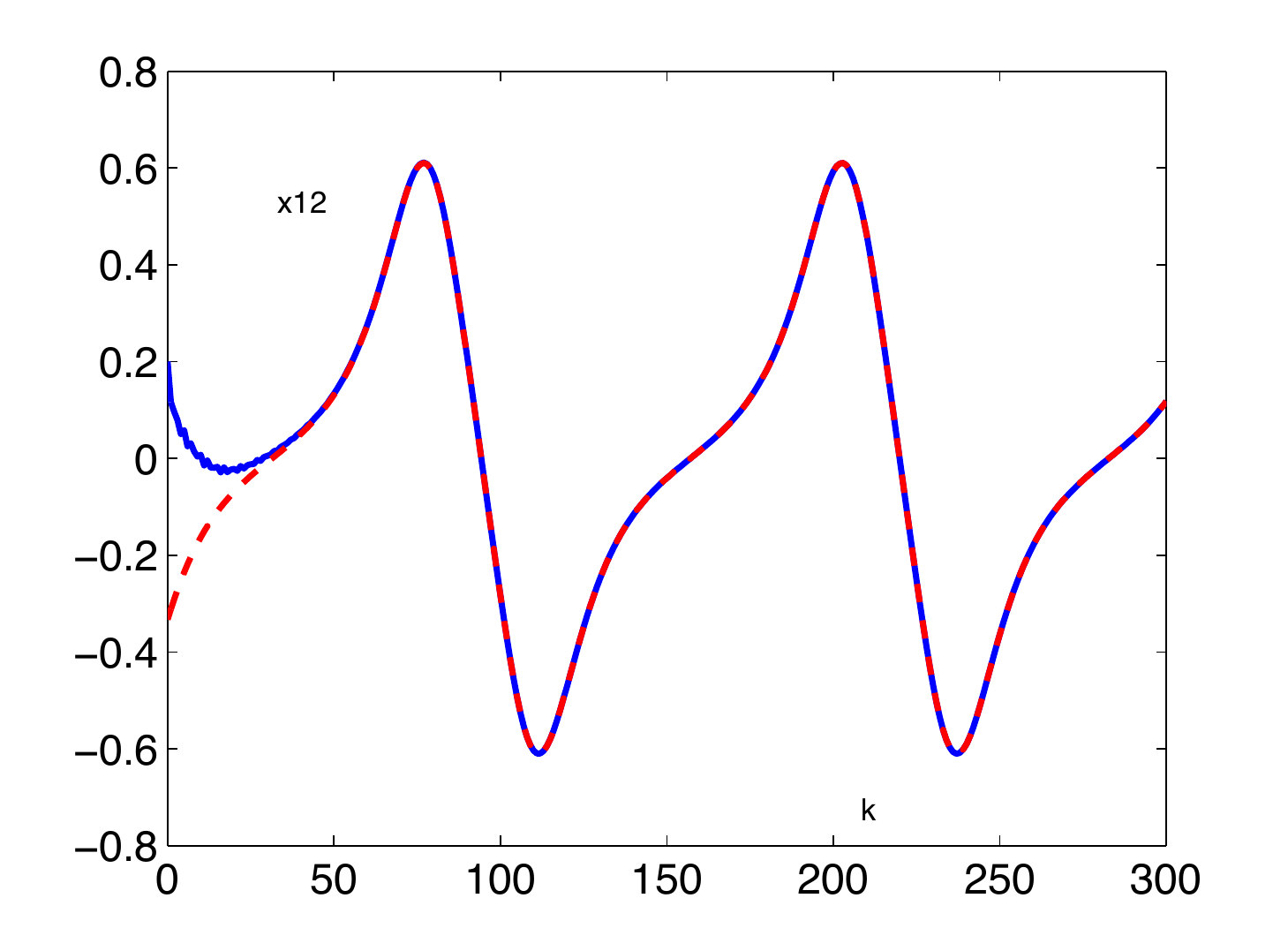}}
\subfigure[Profiles of the (2,1) entry $X_{2,1}$]{\includegraphics[width=0.4\columnwidth]{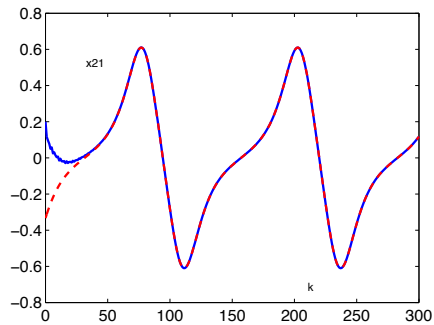}}
\subfigure[Profiles of the (2,2) entry $X_{2,2}$]{\includegraphics[width=0.4\columnwidth]{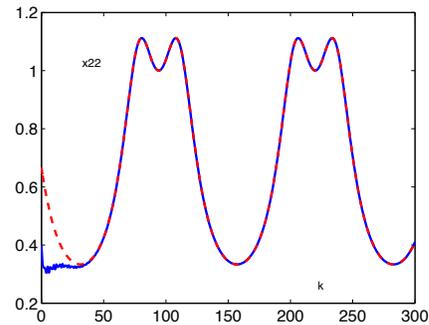}}
\caption{Profiles of the four entries of the solution $X$ when solving the discrete time-varying matrix-inverse problem (\ref{example2}) with $\tau=0.1$ s. Here the solid curves show the solution entries generated by the 5 instants discrete model obtained from random starting values  and the dash-dotted curves depict the theoretical solutions.}
\label{states.i}
\end{figure}
\newpage

\subsection{\large A Zhang Neural Network approach for solving nonlinear convex optimization problems  under time-varying linear constraints} 

This section is based on joint work with  Jian Li, Mingzhi Mao, and Yunong Zhang, see \cite{LMUZh1}.\\[-1mm]

\noindent
{\bf Problem formulation :}\\[-2mm]
\begin{equation*}\label{continuous.problem}
\begin{split}
&\text{find }~~~~~~~~~~~~\text{ min}~~f(x(t),t) \\
&\text{such that }~~~~~A(t)x(t)=b(t) \  \ \ \text{ with } x(t) \in \RR^n, \ b(t) \in \RR^m \text{ and } A(t) \in \RR^{m,n}.
\end{split}
\end{equation*}
{\bf Building a {continuous time model}, i.e.,  formulating an equation for the problem :}\\[1.5mm]
The Zhang Neural Network approach can be built on  the Lagrange function 
 \begin{equation*}
L(x(t),l(t),t)=f(x(t),t)+l^T(t)(A(t)x(t)-b(t)),
\end{equation*}
where $l(t)=[l_1(t),\cdots,l_m(t)]^T\in\RR^m$ is the Lagrange-multiplier vector and  $..^T$ denotes the transpose. Note that there will no need to  {solve for the Lagrange functions  $l(t)$} here. 
Set $$y(t) = [x^T(t), l^T(t)]^T = [y_1(t),\cdots,y_{n}(t),y_{n+1}(t),\cdots,y_{n+m}(t)]^T \in \RR^{n+m}$$ and 
\begin{equation*}
\begin{split}
h(y(t),t)=
\begin{bmatrix}
\dfrac{\partial f(x(t),t)}{\partial x}+A^T(t)l(t)
\\[3mm]
A(t)x(t)-b(t)
\end{bmatrix}=
\begin{bmatrix}
h_1(y(t),t)\\
\vdots\\
h_n(y(t),t)\\
h_{n+1}(y(t),t)\\
\vdots\\
h_{n+m}(y(t),t)
\end{bmatrix}\in \RR^{n+m}.
\end{split}
\end{equation*}
We transform the multiplier problem into an initial value DE problem instead.
By {stipulating exponential decay for $h(t)$} we  obtain the model equation
$$
\dot y(t)=-H^{-1}( y (t),t)\left(\lambda h y(t),t)+\dot h_t(y(t),t)\right)
$$
for the Jacobian matrix
$$
H(y(t),t) =
\begin{bmatrix}
\dfrac{\partial f^2(x(t),t)}{\partial x}\partial ^{T}x &&A^{T}(t)
\\[3mm]
A(t) &&0
\end{bmatrix} \
\text{ and } \ \dot h_t(y(t),t)=
\dfrac{\partial h(y(t),t)}{\partial t}\ .
$$

\vspace*{-1mm}
\noindent
{\bf Discretizing the model and choosing suitable high order {finite difference formulas} :}\\[1mm]
To discretize the continuous model  
\begin{equation*}
\begin{split}
\dot y(t)=-H^{-1}( y (t),t)\left(\lambda h y(t),t)+\dot h_t(y(t),t)\right)
\end{split}
\end{equation*}
we can use the forward {Euler difference formula}  with truncation error order $O(\tau)$
$$ \dot x(t_{k})=\frac{x(t_{k+1})-x(t_{k})}{\tau}$$
 or the {four-instance forward difference formula} (4-IFD)
 $$ \dot x(t_{k})=\frac{5x(t_{k+1})-3x(t_{k})-x(t_{k-1})-x(t_{k-2})}{8\tau} $$
 with truncation error order $O(\tau^2)$. The Euler formula yields the {discretized model}
 $$ y_{k+1}=- H^{-1}(y_k,t_k)\left( \kappa h(y_k,t_k)+\tau\dot h_t(y_k,t_k)\right)+y_{k} \  \text{ with } \ \kappa = \tau \lambda$$
while the 4-IFD formula results in 
 $$ y_{k+1}=-\frac{8}{5} H^{-1}(y_k,t_k)\left( \kappa h(y_k,t_k)+\tau\dot h_t(y_k,t_k)\right)+\frac{3}{5}y_{k}+\frac{1}{5} y_{k-1}+\frac{1}{5}y_{k-2}+O(\tau^3)\ .$$
 Both {discretization formulas} are {consistent and convergent}. 
 This can be proved via the the roots of the associated  characteristic polynomial. Its roots must lie in the complex unit circle and cannot be repeated on its boundary. \\[1mm]
Since the value of $\dot h_t(y_k,t_k)$ may not be known explicitly we may replace it by 
$$ \dot h_t(y_k,t_k)=\frac{3h(y_k,t_k)-4h(y_k,t_{k-1})+h(y_k,t_{k-2})}{2\tau} $$ which uses the three-point backward finite difference formula 
$$\dot x(t_{k})=\frac{3x(t_{k})-4x(t_{k-1})+x(t_{k-2})}{2\tau}\ 
$$
of order $O(\tau^2)$.\\[1mm]
Then the {discretized FIFD formula} becomes more complicated but easier to implement:
\begin{eqnarray*}
 y_{k+1}&=&-\frac{8}{5} H^{-1}(y_k,t_k)\left( \left(\kappa+\frac{3}{2}\right) h(y_k,t_k)-2h(y_k,t_{k-1})+\frac{1}{2}h(y_k,t_{k-2})\right)\\
 && \ \ \ +\frac{3}{5}y_{k}
+\frac{1}{5}y_{k-1}+\frac{1}{5}y_{k-2} \ \quad \text{ of order } \ \ O(\tau^3)\ .
\end{eqnarray*}
To implement this formula, the inverse of the Jacobian matrix $H$ can be computed at each time $t_k$ in a fraction of the available real-time interval $[t_k, t_{k+1})$ by using the real-time inverse finding ZNN method from the previous subsection 2.2.\\[-2mm]

\noindent
{\bf Numerical example and results :}\\[1mm]
As an example we solve  the following convex nonlinear optimization problem with known theoretical solution numerically by using our ZNN method, for further details and applications see \cite{LMUZh1}:\\[-2mm]
\begin{equation*}
\begin{split}
\text{Find } \ \ &\text{min}(\cos(0.1t_{k+1})+2)x_1^2+(\cos(0.1t_{k+1})+2)x_2^2+2\sin(t_{k+1})x_1x_2+\sin(t_{k+1})x_1+\cos(t_{k+1})x_2 \\
&\text{so that}~~~~~\sin(0.2t_{k+1})x_1+\cos(0.2t_{k+1})x_2=\cos(t_{k+1}).
\end{split}
\end{equation*}

\vspace*{-1mm}\enlargethispage{50mm}
\noindent
The FIFD formula is a 4 instance formula, while the Euler formula needs only two. 
Both discretization models work in real-time and both typically create the optimal solution  in a fraction of a second with differing degrees of accuracy according to their orders.\\[1mm]
The example below runs for 10 sec. The time-varying values for $f(x(t),t)$, $A(t)$ and $b(t)$ are given as functions  and evaluated from their function formulations. 
 In real world applications these values might be supplied by sensors during each time interval $t_i \leq t_{i+1}$ and the empirical  values would be inserted into the difference formulas as they are evaluated by sensors in real time. 
 
\vspace*{-5mm}
\begin{figure}[H]\centering
\psfrag{t}[c][c][1] {$k\tau$}%
\psfrag{eee}[c][c][1]{$e(k\tau)$}%
\psfrag{x1}[c][c][1]{$x_1$}%
\psfrag{x2}[c][c][1]{$x_2$}%
\psfrag{NFIFD}[c][c][0.8]{~~~~~~~FDZTND-U}%
\psfrag{Euler}[c][c][0.8]{~~~~~~~~~EDZTND-U}%
\psfrag{Theory}[c][c][1]{~~~~~$\mathbf{x}^*(k\tau)$}%
\subfigure[]{\includegraphics[width=0.4\columnwidth]{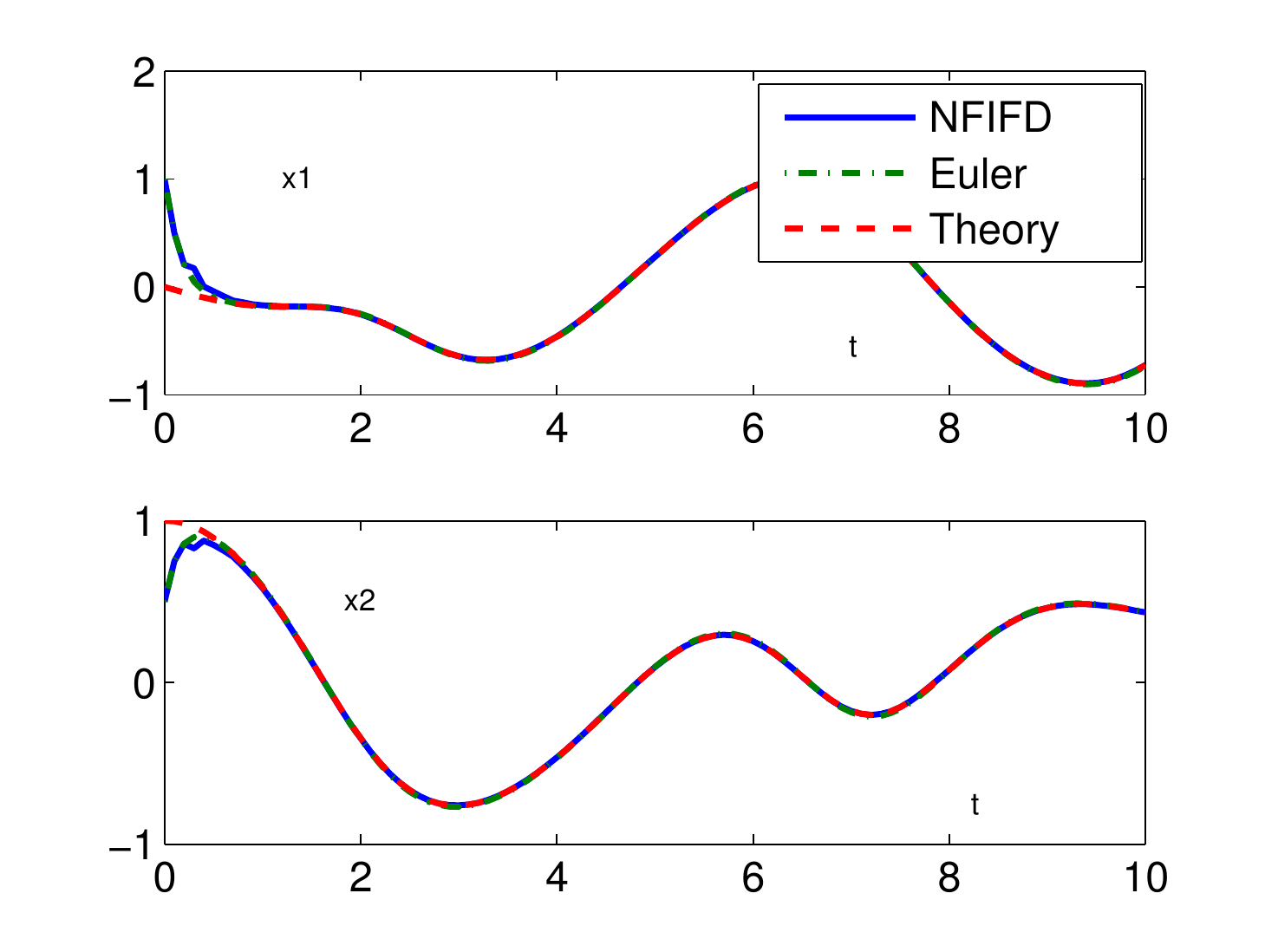}\label{fig.dt.0.1.states.u}}
\subfigure[]{\includegraphics[width=0.4\columnwidth]{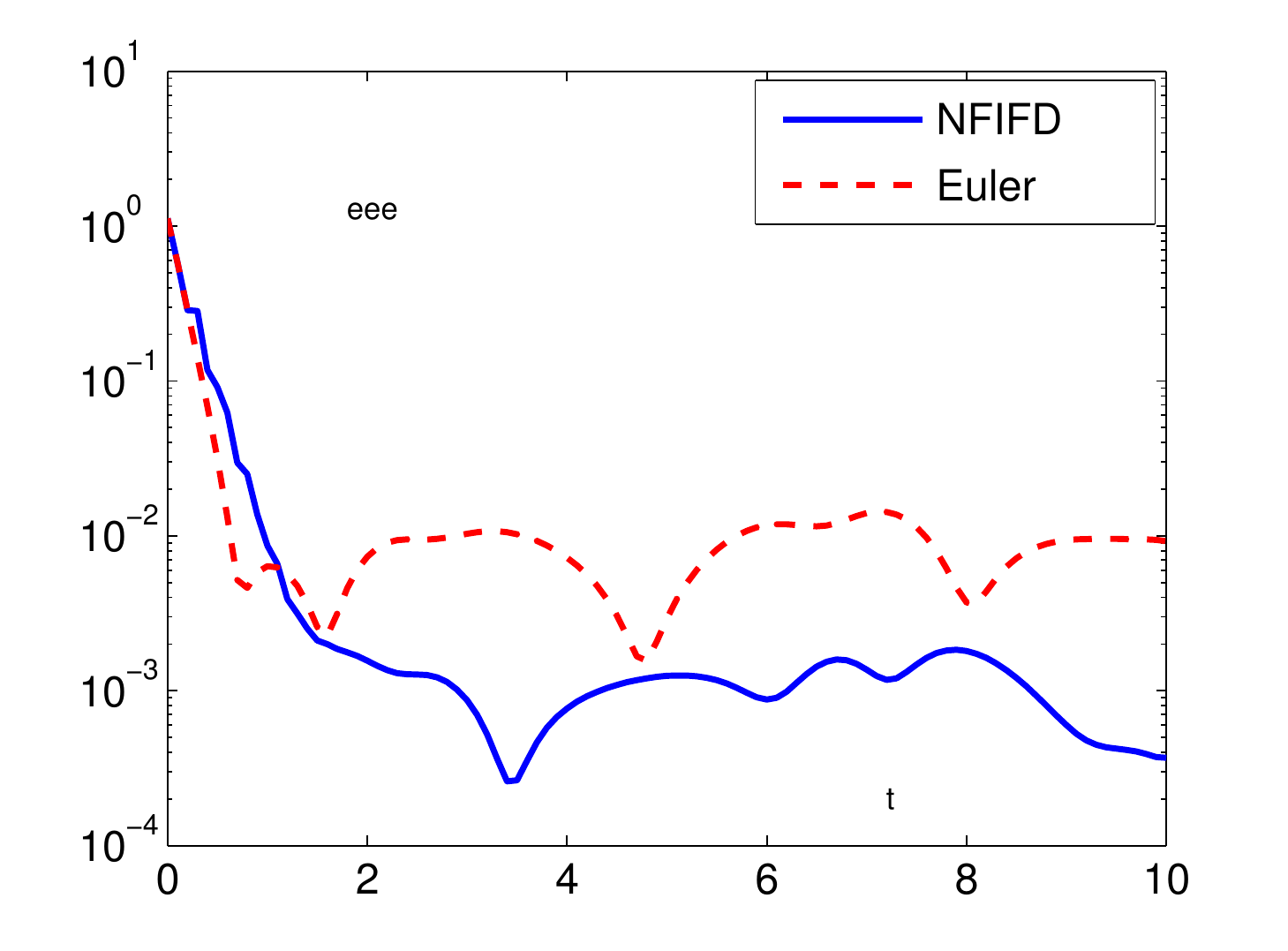}\label{fig.dt.0.1.lambda.u}}
\subfigure[]{\includegraphics[width=0.4\columnwidth]{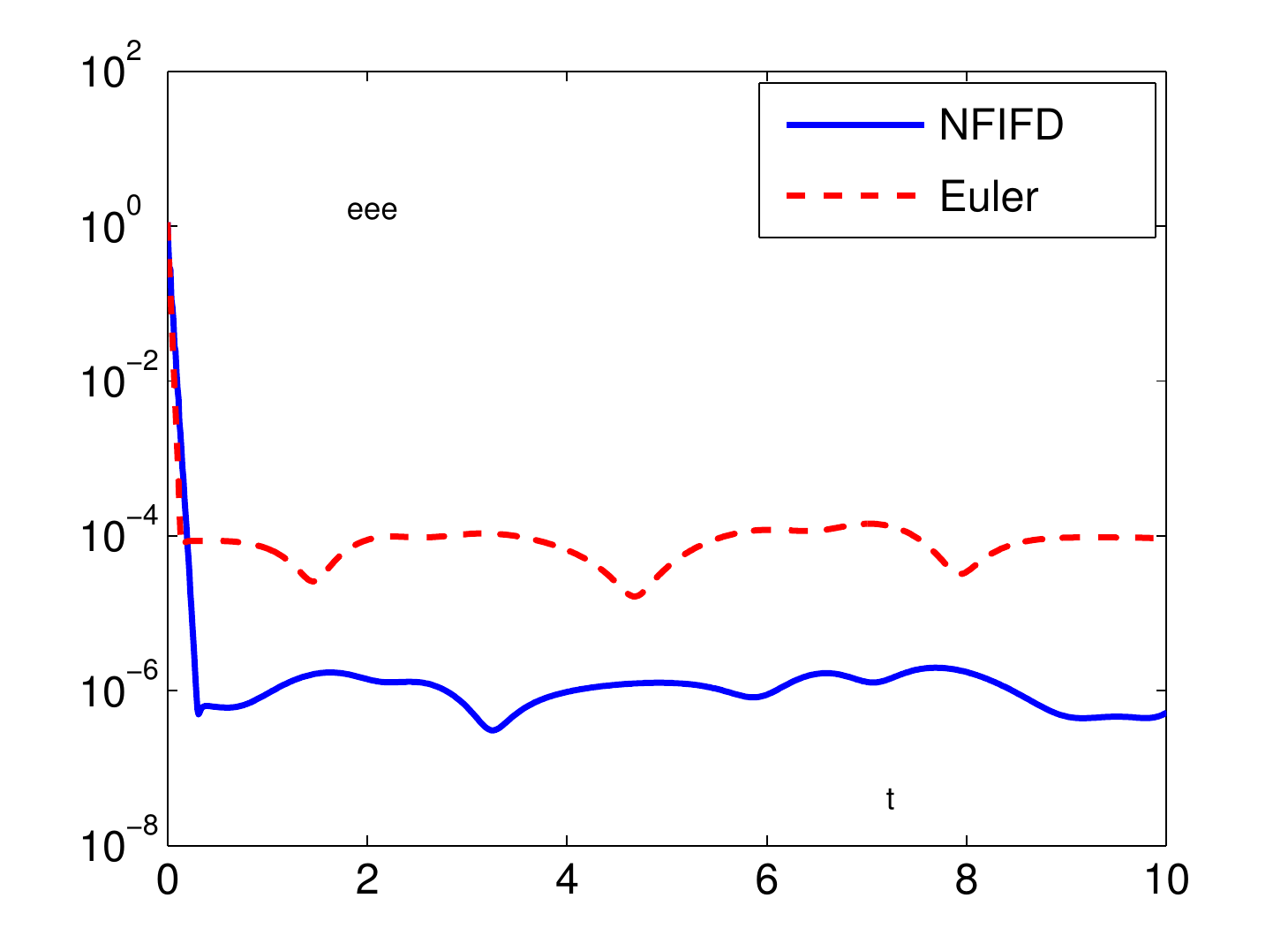}\label{fig.dt.0.01.lambda.u}}
\subfigure[]{\includegraphics[width=0.4\columnwidth]{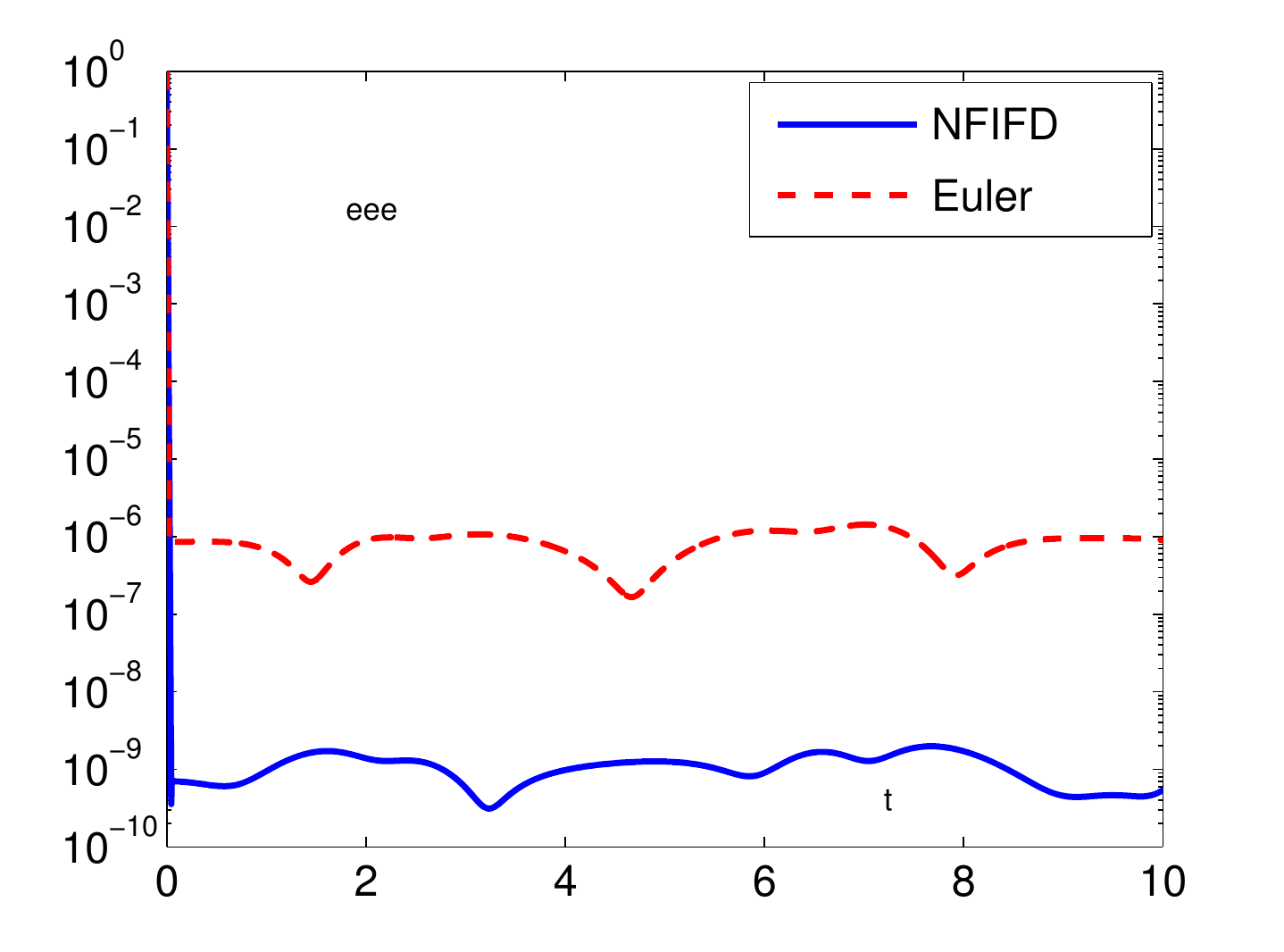}\label{fig.dt.0.001.lambda.u}}
\caption{ {Solution states} with $\tau=0.1$ sec, and {solution errors} generated by the 4-IFD based discretization model  (of order  $O(\tau^3)$) and the Euler based model (of order $O(\tau^2)$) with $\tau = 0.1$,  $0.01$, and $0.001$ sec and $\lambda = 10$:\newline
\hspace*{26mm} (a) {solution states} with $\tau=0.1$ sec, \qquad (b) {solution errors} with $\tau=0.1$ sec, \newline
\hspace*{26mm} (c) {solution errors} with $\tau=0.01$ sec,\hspace*{4.76mm}
 (d) {solution errors} with $\tau=0.001$ sec.}
\end{figure}
\newpage

\vspace*{-1mm}
\section{\Large On Quantum and Multi-state Computing \ \  (\emph{Epoch 8} \ \ \ \normalsize{yet to start and come}\Large{)}}

\vspace*{-1mm}
{Quantum Computing} and  {multi-state memory and computers} with multi-state processors will change the way we compute once they become available. 
 They will require new operating systems and new software with new and yet to be discovered algorithms. 
What will this new era entail? Nobody knows or can reliably predict.\\[1mm]
 I asked an `expert' on quantum computing three years ago as to when he expected to have a quantum computer at his disposal or on his desk. 
 The answer was : ``Not in my lifetime, not in 20 years".\\[1mm]
 Currently about a dozen or more research centers in Europe and South-East Asia are trying to build quantum computers based on the quantum superposition principle and quantum entanglement of elementary particles. They do so  in a multitude of different ways.
The envisioned benefit of these efforts would be to be able to compute super fast in parallel and in simulations to solve huge data problems quicker than ever before and to solve problems  that are unassailable now with our current best supercomputer networks. All of the proposed  quantum science techniques make use of superconducting circuits and particles. The aim is to build quantum computers in one or two decades with around 100 entangled quantum bits. Such a quantum computer would be bulky, it would need much supplementary  equipment for cooling and so forth and could easily take up a whole floor of a building,  just as the first German and British valve computers did in the 1940s. But it would surpass the computing capacity of all current supercomputers and desk and laptops on Earth combined. Currently the largest working entangled quantum array contains fewer than 10 quantum bits. Access of  a 100 bit quantum computer would probably be via the Cloud and there would be no quantum computer laptops.
Quantum computers may take another 10, 20 or 30 years to materialize. \\
 \emph{How will they come about? Which yet unknown algorithms will they use? Who will invent them? Who code them?} \\
If history can be a guide,  John Francis and Vera Kublanovskaya were both working  independently on circuit diagrams and logic gate designs for valve computers in England and in Russia at the time when they discovered QR (or LQ) in the late 1950s. \\[1mm]
So we possibly are looking for quantum computer hardware and software designers who know numerical analysis and algorithm development in or about the year 2040. 
In a similar fashion, Leibnitz and Seki formalized our now ubiquitous matrix concept independently but simultaneously in   1683, in Germany and in Japan.\\[2mm]
\hspace*{2mm} \hfill \emph{Maybe it will take two again?}  \\[4mm]

\noindent
The references given below only go back to the year 1799.\\[-8mm]

 \vspace*{10mm}
 
\bigskip

\bigskip

\end{document}